\documentclass{amsart}
\usepackage{amssymb,amsmath}

\newtheorem{theorem}{Theorem}[section]
\newtheorem{lemma}[theorem]{Lemma}

\newtheorem{cor}[theorem]{Corollary}
\newtheorem*{unthm}{Theorem}

\theoremstyle{definition}

\theoremstyle{remark}
\newtheorem{remark}[theorem]{Remark}

\def\a{{\mathfrak a}}

\def\g{{\mathfrak g}}
\def\h{{\mathfrak h}}
\def\k{{\mathfrak k}}
\def\l{{\mathfrak l}}
\def\n{{\mathfrak n}}

\def\p{{\mathfrak p}}
\def\q{{\mathfrak q}}
\def\r{{\mathfrak r}}

\def\t{{\mathfrak t}}

\def\z{{\mathfrak z}}

\def\C{{\mathbb C}}
\def\R{{\mathbb R}}

\numberwithin{equation}{section}

\begin{document}
\author[M. Otto]{Michael Otto}
\address{Department of Mathematics, 
University of Arizona, Tucson, AZ 85721.}
\email{}
\title[Restriction of the moment map]{Restriction of the moment map to certain non-Lagrangian submanifolds}

\begin{abstract} 
Let $T\times M \rightarrow M$ be a Hamiltonian torus action on a connected symplectic manifold $M$ for which the associated moment map $\Phi:M\rightarrow \t^*$ is proper as a map into a convex open set $\rho \subseteq \t^*$. 
We consider a closed submanifold $Q$ of $M$ and show that under certain local conditions on $Q$ one has $\Phi(Q)=\Phi(M)$. We apply this result in the special case that $Q$ arises as the fixed point set of some involution $\sigma$ on $M$ which is not necessarily antisymplectic. 
\end {abstract}
\maketitle

\section{Introduction}

Let $\Phi:M\rightarrow \t^*$ be the moment map for a Hamiltonian torus action $T\times M \rightarrow M$ on a connected symplectic manifold $(M,\omega)$. The classical symplectic convexity theorem of Atiyah-Guillemin-Sternberg \cite{A,GS} asserts that for compact $M$ the image $\Phi(M)$ is a convex set. More precisely, $\Phi(M)=conv(\Phi(M^T))$, i.e. the moment map image is the convex hull of the (finitely many) images of $T$-fixed points of $M$. Duistermaat \cite{D} added an involutive diffeomorphism $\sigma:M\rightarrow M$ to the picture. 

\begin{unthm}(Duistermaat) 
Let $M$ be compact. In addition, let $\sigma:M\rightarrow M$ be an antisymplectic involution which leaves the moment map invariant, i.e. $\sigma^*\omega=-\omega$ and $\Phi \circ \sigma = \Phi$. Suppose that the $\sigma$-fixed set $Q=M^\sigma$ is not empty. \\
Then, $\Phi(Q_0)=\Phi(M)$ for any connected component $Q_0$ of $Q$. In particular, $\Phi(Q)=\Phi(M)$.
\end{unthm}

The AGS-theorem and Duistermaat's extension have been used, among other things, to provide a symplectic proof of Kostant's nonlinear convexity theorem for semisimple Lie groups (see \cite{LR}). However, the condition that $\sigma$ be antisymplectic in Duistermaat's theorem turned out to be too restrictive to deal with all instances of Kostant's theorem. This gave the motivation for two different generalizations of Duistermaat's result. In his thesis \cite{Sl}, Sleewaegen took a very general approach and considered a closed submanifold $Q$ of $M$, independent of any involution $\sigma$. He found conditions on the restriction of the moment map $\Phi|_Q$ sufficient to ensure $\Phi(Q)=\Phi(M)$. These conditions involve gradients and Hessians of Hamiltonian functions $\Phi_X$ and their restrictions $\Phi_X|_Q$. In \cite{KO}, on the other hand, the authors replaced the assumption $\sigma^* \omega=-\omega $ in Duistermaat's theorem with a weaker condition and showed that $\Phi(Q)=\Phi(M)$ still holds true.

This paper is aimed at providing a symplectic convexity theorem for submanifolds $Q$ of $M$ which generalizes the results in \cite{Sl} and \cite{KO} (and Duistermaat's theorem). Most importantly, $M$ is no longer required to be compact. Instead we assume that the moment map is proper in the sense of \cite{LMTW}. This means that $\Phi:M\rightarrow \t^*$ is proper as a map into an open convex set $\rho \subseteq \t^*$. In this case, $\Phi(M)$ still is a convex set as shown in \cite{LMTW}. We now state the main result. 

\begin{unthm} 
Let $T\times M \rightarrow M$ be a Hamiltonian torus action on a connected symplectic manifold $M$. Suppose the associated moment map $\Phi:M\rightarrow \t^*$ is proper as a map into a convex open set $\rho \subseteq \t^*$. Let $Q$ be a nonempty closed submanifold of $M$ such that the following two conditions are satisfied.
\begin{align*}
(A) & \quad \text{The derivative $D\Phi|_Q(q)$ is surjective at each point $q\in Q$ at which} \\
& \quad \text{$D\Phi(q)$ is surjective.} \\
(B) & \quad \text{Let $q\in Q, \ X\in \t_q$ (the stabilizer Lie algebra of $q$). If for each neighborhood} \\
& \quad \text{$U\subseteq M$ of $q$ there exists some $p\in U$ with $\Phi_X(p)>\Phi_X(q)$, then each such $U$} \\
& \quad \text{contains some $q'\in Q$ with $\Phi_X(q')>\Phi_X(q)$.} 
\end{align*} 
Then, $\ \Phi(Q)=\Phi(M)$.
\end{unthm}

The theorem is proved in Section \ref{main}. Throughout the proof we make extensive use of the local normal form of the moment map. We also exploit the properness of $\Phi$ and study its regular values. In Section \ref{invol} we apply the theorem to the situation where $Q=M^\sigma$ is the fixed point set of some involution $\sigma:M\rightarrow M$. We obtain $\Phi(Q)=\Phi(M)$ under rather general conditions on $\sigma$. In particular, it is not necessary that $Q$ is Lagrangian. Finally, in Section \ref{subvdb} we use the results of \ref{invol} to give a symplectic proof of van den Ban's convexity theorem for semisimple symmetric spaces. This is a generalization of Kostant's theorem involving noncompact group orbits.

\section{The main result}\label{main}

Consider a Hamiltonian torus action $T\times M \rightarrow M$ on a connected symplectic manifold $M$. We assume that the associated moment map $\Phi:M\rightarrow \t^*$ is proper as a map into a convex open set $\rho \subseteq \t^*$. In this case the moment map image $\Phi(M)$ is a convex set (see \cite{LMTW}). Without loss of generalization we assume that the interior of $\Phi(M)$ in $\t^*$ is nonempty (this can always be achieved by factoring out a suitable subtorus of $T$). 

Let $Q$ be a (nonempty) closed submanifold of $M$. The goal is to establish conditions on $Q$ so that its moment map image $\Phi(Q)$ is all of $\Phi(M)$. First we define regular elements $M_{reg}=\{ p\in M : \t_p = 0 \}$, where $\t_p$ denotes the Lie algebra of the stabilizer $T_p\subseteq T$ of $p$. We set $Q_{reg}=Q\cap M_{reg}$. For later use we also define the set of regular and critical values of $\Phi:M\rightarrow \rho$ in the usual way.
\[ V_{crit}=\{ X\in \rho : \exists \ p\in \Phi^{-1}(p) \ \text{such that} \ D\Phi(p) \ \text{is not surjective} \} , \]
and we set $V_{reg}=\Phi(M)\setminus V_{crit}$.
Recall that each $X\in \t$ defines the Hamiltonian function $\Phi_X:M\rightarrow \R $ and its restriction $\Phi_X|_Q:Q\rightarrow \R $. A point $p\in M$ is a critical point of $\Phi_X$ if and only if $X\in \t_p$. Also $p\in M$ is regular if and only if $D\Phi(p)$ is surjective.

Consider the following two conditions on $Q$. 
\begin{align*}
(A) & \quad \text{The derivative $D\Phi|_Q(q)$ is surjective at each $q\in Q_{reg}$.} \\
(B) & \quad \text{Let $q\in Q, \ X\in \t_q$. If for each neighborhood $U\subseteq M$ of $q$ there exists} \\
& \quad \text{some $p\in U$ with $\Phi_X(p)>\Phi_X(q)$, then each such $U$ contains some} \\
& \quad \text{$q'\in Q$ with $\Phi_X(q')>\Phi_X(q)$.}
\end{align*} 

\begin{remark}
Compare (A) and (B) with Def. and Prop. 3.4.1 in \cite{Sl}. The condition on the Hessians of $\Phi_X$ and $\Phi_X|_Q$ in \cite{Sl} is replaced with the weaker "topological" condition (B).
\end{remark}

The goal of this section is to show that conditions (A) and (B) are sufficient to ensure $\Phi(Q)=\Phi(M)$. 

First we study local properties of $\Phi(M)$ and $\Phi(Q)$. The main ingredient will be a local normal form for the moment map. Fix $p\in M$ with stabilizer Lie algebra $\t_p$. Consider the decomposition $\t^*=\t_p^\perp \oplus \t_p^*$. Recall the local normal form description of the moment map $\Phi$ (see e.g. \cite{HNP}, Lemma 2.2). In appropriate (symplectic) coordinates $(v,w,x,y)\in \R^m \times \R^m \times \R^n \times \R^n$ on a neighborhood $U_p$ of $p=(0,0,0,0)$, the moment map takes on the form
\begin{equation}\label{nf}
\Phi(v,w,x,y)=(v,\frac{1}{2}\sum_{i=1}^{n} \alpha_i (x_i^2+y_i^2)) \in \R^m \times \t_p^* .
\end{equation}
Here, $m=dim(\t_p^\perp)$ and $2m+2n=dim(M)$. Without loss of generality we have assumed that $\Phi(p)=0$. The $\alpha_i\in \t_p^*$ are the weights of the linearized action of the stabilizer $T_p$.

The normal form (\ref{nf}) is useful to describe $\Phi$ on a neighborhood $U_p$ of $p\in M$. However, to study $\Phi(M)$ on a neighborhood of $\Phi(p)$ we need to combine the local information on neighborhoods $U_m$ for several $m\in \Phi^{-1}(\Phi(p))$. The properness of $\Phi$ ensures that we need to consider only finitely many such points $m$, as the following lemma shows.

\begin{lemma}\label{finnbhds}
Fix $p\in M$. There exist a neighborhood $U$ of $\Phi(p)$ in $\rho$ and finitely many points $p_1, \dots , p_k \in \Phi^{-1}(\Phi(p))$ with neighborhoods $U_1, \dots , U_k \subseteq M$, respectively, with the following properties. 
\begin{enumerate}
\item For each $i$, the restriction $\Phi|_{U_i}$ takes on a local normal form as in (\ref{nf}) where $U_i$ is centered at $p_i$, i.e. $p_i=(0,0,0,0)$ in local coordinates. \\
\item $\Phi^{-1}(U)\subseteq U_1 \cup \dots \cup U_k$.
\end{enumerate}
\end{lemma}
\begin{proof}
Fix a compact neighborhood $V$ of $\Phi(p)$. Since $\Phi$ is proper $\Phi^{-1}(V)$ is compact as well. For each $x\in \Phi^{-1}(V)$ we choose an open neighborhood $U_x$ small enough so that $\Phi|_{U_x}$ has a local normal form as in (\ref{nf}). If $x\not\in \Phi^{-1}(p)$ we shrink $U_x$, if necessary, to ensure that $\Phi(U_x)$ is contained in some open set $V_x\subseteq \rho$, for which $\Phi(p)\not\in \overline{V_x}$. 
The compact set $\Phi^{-1}(V)$ can now be covered by finitely many such open sets $U_x$. We write $U_1, \dots , U_k$ for those that are centered at points $p_1, \dots , p_k \in \Phi^{-1}(\Phi(p))$, and $U_{p_k+1}, \dots , U_{p_l}$ for those centered at $p_{k+1}, \dots , p_l \not\in \Phi^{-1}(\Phi(p))$. Now, choose a neighborhood $U\subseteq V$ of $\Phi(p)$ such that $U\cap \overline{V_{p_i}}=\emptyset$ for all $i=k+1, \dots , l$. 
For $m\in \Phi^{-1}(U)$, suppose $m\in U_{p_i}$ for some $i\in \{ k+1, \dots , l \} $. Then $\Phi(m)\in \Phi(U_{p_i})\subseteq V_{p_i}$ and $\Phi(m)\in U$, a contradiction. This shows that $\Phi^{-1}(U)\subseteq U_1 \cup \dots U_k$. 
\end{proof}

With Lemma \ref{finnbhds} and the local normal form (\ref{nf}) at hand we can describe the critical values near a point $\Phi(p)$. Choose neighborhoods $U_1, \dots , U_k$ centered at $p_1,\dots,p_k \in \Phi^{-1}(\Phi(p))$ and $U\subseteq \rho$ as in Lemma \ref{finnbhds}. On each $U_i$ one defines the local moment cone $\Gamma_i=\t_{p_i}^\perp + cone\{\alpha_1^i,\dots,\alpha_{n_i}^i\}\subseteq \t^*$, where $cone\{\alpha_1^i,\dots,\alpha_{n_i}^i\}$ denotes the cone spanned by the weights $\alpha_1^i,\dots,\alpha_{n_i}^i\in \t_{p_i}^*$ that occur in the normal form (\ref{nf}) of $\Phi|_{U_i}$. For a sufficiently small neighborhood $U'$ of $\Phi(p)$ one has $\Phi(U_i)\cap U'=(\Phi(p)+\Gamma_i)\cap U'$. The important fact is that the local moment cone $\Gamma_i$ is the same for all $i=1,\dots ,k$. To see this note that a local moment cone $\Gamma_m$ can be defined at each $m\in \Phi^{-1}(\Phi(p))$. From the normal form for $\Phi|_{U_i}$ one can see that $\Gamma_m$ is the same for all $m\in \Phi^{-1}(\Phi(p))\cap U_i$ (even though stabilizer Lie algebra $\t_m$ and weights might differ). In addition, the fiber $\Phi^{-1}(\Phi(p))$ is connected (this was also shown in \cite{LMTW}) and covered by $U_1\cup \dots \cup U_k$. Therefore, $\Gamma_1=\dots =\Gamma_k=:\Gamma$. Due to our general assumption $int(\Phi(M))\neq\emptyset$ we know $\t_{p_i}^*=span\{\alpha_1^i,\dots,\alpha_{n_i}^i\}$, i.e. $int(\Gamma)\subseteq \t^*$ is non-empty as well. 
Fix $i$ for a moment. Using the weights $\alpha_1^i,\dots,\alpha_{n_i}^i\in \t_{p_i}^*$ we define certain cones of codimension one in $\t^*$ in the following way. Let $\{\alpha_{j_1}^i,\dots,\alpha_{j_r}^i\}$ be a maximal subset of weights such that the subspace $H=\t_{p_i}^\perp \oplus span\{\alpha_{j_1}^i,\dots,\alpha_{j_r}^i\}$ has codimension one in $\t^*$. We define the cone $\gamma=\t_{p_i}^\perp + cone\{\alpha_{j_1}^i,\dots,\alpha_{j_r}^i\}$. Enumerate all cones $\gamma_1^i,\dots,\gamma_{s_i}^i$ that can be obtained this way from the set of weights $\{\alpha_1^i,\dots,\alpha_{n_i}^i\}$. Note that if $dim(\t_{p_i})=1$ there is only one such cone $\gamma_1^i=\t_{p_i}^\perp$. If $dim(\t_{p_i})=0$ there is no such cone.

\begin{lemma}\label{ucrit}
For $p\in M$, let $U,U_1,\dots,U_k,p_1,\dots,p_k$ be as in Lemma \ref{finnbhds}. Without loss of generalization $U$ is chosen sufficiently small such that $U\subseteq \Phi(U_i)$ for all $i=1,\dots,k$. Let $\gamma_1^1,\dots,\gamma_{s_1}^1,\dots,\gamma_1^k,\dots,\gamma_{s_k}^k$ be the cones of codimension one defined in the paragraph preceding this lemma. 
\begin{enumerate}
\item Then, $x\in U$ is a critical value if and only if $x\in \bigcup_{i=1}^k \bigcup_{j=1}^{s_i} \gamma_j^i$. \\
\item Let $\Gamma$ be the local moment cone at $\Phi(p)$. Fix a ball $B\subseteq U$ centered at $\Phi(p)$. \\
Then, $\Phi^{-1}(B)\subseteq U_1\cup \dots \cup U_k$. Moreover, the set $B\cap(\Phi(p)+\Gamma)=B\cap \Phi(M)$ is partitioned by $\{ \gamma^i_j \} $ into open (and convex) chambers, say $C_1,\dots ,C_N$. These chambers are exactly the connected components of $V_{reg}\cap B\cap (\Phi(p)+\Gamma)=V_{reg}\cap B$.
\end{enumerate}
\end{lemma}
\begin{proof}
\begin{enumerate}
\item Suppose $x\in \gamma_j^i$ for some $i,j$. Then $x\in \t_{p_i}^\perp + cone\{\alpha_{j_1}^i,\dots,\alpha_{j_r}^i\}$ for some maximal subset $\{\alpha_{j_1}^i,\dots,\alpha_{j_r}^i\}\subseteq \t_{p_i}^*$ for which $span\{\alpha_{j_1}^i,\dots,\alpha_{j_r}^i\}$ has codimension one in $\t_{p_i}^*$. Fix a nonzero $X\in \t_{p_i}\cap(span\{\alpha_{j_1}^i,\dots,\alpha_{j_r}^i\})^\perp$. Consider a preimage point $m\in U_i\cap \Phi^{-1}(x)$. The local normal form (\ref{nf}) for $\Phi$ on $U_i$ shows that the entries $(x_d,y_d)$ of $m=(v,w,x,y)$ may be nonzero only for $d\in \{j_1,\dots,j_r\}$, since $x\in \gamma_j^i$. But then $m$ is fixed under the action of the one-parameter-group $\{ \exp tX :t\in \R \}$. This implies $D\Phi_X(m)=0$, i.e. $D\Phi(m)$ is not surjective. This shows that $x=\Phi(m)$ is a critical value. 

Conversely, if $x$ is critical there must be a preimage point $m$ with $D\Phi(m)$ not surjective. Fix a nonzero $X\in \t$ with $D\Phi_X(m)=0$, i.e. $m$ is fixed under $\{ \exp tX :t\in \R \}$. According to Lemma \ref{finnbhds} the preimage point $m$ is contained in some $U_i$. Suppose that in local coordinates $m=(v,w,x,y)$ the pairs $(x_d,y_d)$ that are nonzero correspond to the weights $\{\alpha_{j_1}^i,\dots,\alpha_{j_r}^i\}$. Since $m$ is fixed under $X$ we must have $\alpha_{j_1}^i(X)=\dots=\alpha_{j_r}^i(X)=0$. Therefore, $span\{\alpha_{j_1}^i,\dots,\alpha_{j_r}^i\}\subsetneq \t_{p_i}^*$, and $\{\alpha_{j_1}^i,\dots,\alpha_{j_r}^i\}$ is a subset of a maximal set of weights that defines a cone $\gamma_j^i$. But then $x=\Phi(m)\in \t_{p_i}^\perp + cone\{\alpha_{j_1}^i,\dots,\alpha_{j_r}^i\} \subseteq \gamma_j^i $. \\
\item The statement about the chambers is just a reformulation of part (1). The inclusion $\Phi^{-1}(B)\subseteq U_1\cup \dots \cup U_k$ follows directly from property (2) in Lemma \ref{finnbhds}.
\end{enumerate}
\end{proof}

We now focus on the nonempty closed submanifold $Q\subseteq M$. Recall conditions (A) and (B).

\begin{lemma}\label{dense}
Suppose $Q$ satisfies (B). Then, the following assertions hold.
\begin{enumerate}
\item Let $q\in Q$ with $dim(\t_q)>0$ and $U$ a neighborhood of $q$ in $M$. Then, $U$ contains some $q'\in Q$ with $dim(\t_{q'})< dim(\t_q)$. \\
\item $Q_{reg}$ is dense in $Q$. 
\end{enumerate}
\end{lemma}
\begin{proof}
\begin{enumerate}
\item Recall the local normal form of the moment map on a neighborhood $U_q\subseteq U$ of $q\in M$ as in (\ref{nf}). As mentioned earlier, due to our assumption $int(\Phi(M))\neq\emptyset$, we know $\t_q^*=span\{\alpha_1,\dots,\alpha_n\}$. 
In particular, there exists $X\in \t_q$ such that $\Phi_X(U_q)$ contains an open interval $(-\epsilon,\epsilon)$. Therefore, each neighborhood $U$ of $q$ contains some point $p\in M$ with $\Phi_X(p)>\Phi_X(q)=0$. Condition (B) assures the existence of some $q'\in Q\cap U_q$ with $\Phi_X(q')> 0$. From the local normal form it is clear that such a $q'$ is not fixed under the one parameter group $\{ \exp(tX):t\in \R \} $. On the other hand, $\t_p \subseteq \t_q$ for all $p\in U_q$. Therefore, $\t_{q'}\subset \t_q$. \\
\item This is a consequence of part (1).
\end{enumerate}
\end{proof} 

The next lemma is central.

\begin{lemma}\label{inner}
Suppose $Q$ satisfies (A) and (B). Let $q\in Q$ be such that $\Phi(q)$ is an inner point of $\Phi(M)\subseteq\rho$. \\
Then $\Phi(q)$ is an inner point of $\Phi(Q)$.
\end{lemma}
\begin{proof}
If $dim(\t_q)=0$, then $q\in Q_{reg}$, and property (A) implies that $\Phi(q)\in int(\Phi(Q))$. 

If $dim(\t_q)>0$, we consider an open ball $B$ around $\Phi(q)$ as in part (2) of Lemma \ref{ucrit}. The local moment cone $\Gamma$ is all of $\t^*$, since $\Phi(q)\in int(\Phi(M))$ by assumption. The cones $\{ \gamma_j^i \}$ partition $B$ into finitely many open and convex chambers $C_1,\dots,C_N$, all of which are connected components of $V_{reg}\cap B$.
Part (2) of Lemma \ref{dense} shows that there exists some $q'\in Q_{reg}\cap\Phi^{-1}(B)$. Since $C_1\cup \dots\cup C_N$ is dense in $B$, the openness of $Q_{reg}$ now implies that there is some $q''\in Q$ with $\Phi(q'')\in C_a$ for some $a\in \{1,\dots,N\}$. The nonempty set $\Phi(Q)\cap C_a$ is closed in $C_a$, since $\Phi$ is proper. To see that $\Phi(Q)\cap C_a$ is open in $C_a$ fix some $x\in \Phi(Q)\cap C_a$. If $x=\Phi(\tilde{q})$ for some $\tilde{q}\in Q$ the fact that $x\in C_a$ is a regular value implies $\tilde{q}\in Q_{reg}$. Again, the openness condition (A) shows that $\Phi(Q)$ contains an open neighborhood of $\Phi(\tilde{q})$. Hence, $\Phi(Q)\cap C_a$ is open since $C_a$ is. As $C_a$ is connected as well we conclude $\Phi(Q)\cap C_a=C_a$.
With the same argument we can deduce for each chamber $C_b$: If $\Phi(Q)\cap C_b\neq \emptyset$ then $\Phi(Q)\cap C_b=C_b$. Therefore, to show $C_1\cup \dots \cup C_N\subseteq \Phi(Q)$ it is sufficient to show that $\Phi(Q)\cap C_b \neq \emptyset \ \forall b$.
Suppose that, rearranging indices if necessary, $\Phi(Q)\cap C_b\neq \emptyset$ for $b=1,\dots,K$, and $\Phi(Q)\cap C_b =\emptyset$ for $b=K+1,\dots,N$, and that the chambers $C_1$ and $C_N$ have a common wall $\gamma_j^i$. Fix a point $x\in \gamma_j^i$ such that a neighborhood $U_x$ is contained in $\overline{C_1}\cup\overline{C_N}$, but $U_x\cap\overline{C_b}=\emptyset$ for all $b\in\{2,\dots,N-1\}$. Since $\Phi(Q)$ is closed, $\overline{C_1}\subseteq \Phi(Q)$, so there is some $q_0\in Q$ with $\Phi(q_0)=x$. The inclusion $\Phi^{-1}(B)\subseteq U_1 \cup \dots \cup U_k$ implies that $q_0\in Q\cap U_l$ for some $l\in \{ 1,\dots,k \} $. There are two possibilities.

Case 1: The common wall $\gamma_j^i$ can be obtained from weights in the normal from of $\Phi|_{U_l}$ for at least one such $l$. Then $\gamma_j^i=\t_{q_0}^\perp +cone\{\alpha_{j_1}^l,\dots,\alpha_{j_r}^l\}$, where $\{\alpha_{j_1}^l,\dots,\alpha_{j_r}^l\}$ is a maximal subset of the weights $\{\alpha_1^l,\dots,\alpha_{n_l}^l\}$ occuring in the normal form of $\Phi|_{U_l}$ such that $\t_{q_0}^\perp +span\{\alpha_{j_1}^l,\dots,\alpha_{j_r}^l\}\subseteq \t^*$ has codimension one. Fix some nonzero $X\in \t$ that is perpendicular to $\gamma_j^i$. From the normal form for $\Phi|_{U_l}$ we see that $q_0$ is a critical point for $\Phi_X$, i.e. $X\in \t_{q_0}$. Since $x=\Phi(q_0)\in B$ is an inner point of $\Phi(M)$ the image $\Phi_X(U_l)\in \R$ contains an open interval $(-\epsilon,\epsilon)$. Hence, each neighborhood $U$ of $q$ contains points $p',p''\in U\cap U_l$ such that $\Phi_X(p')>\Phi_X(q_0)$ and $\Phi_{-X}(p'')>\Phi_{-X}(q_0)$. Condition (B) asserts the existence of $q_0',q_0''\in Q\cap U_l$ with $\Phi_X(q_0')>\Phi_X(q_0)>\Phi_X(q_0'')$. Recall that $\gamma_j^i$ is part of a hyperplane separating $C_1$ and $C_N$ and that $X$ is perpendicular to $\gamma_j^i$. Therefore, either $\Phi(q_0')$ or $\Phi(q_0'')$ must lie in $C_N$. This contradicts the assumption $\Phi(Q)\cap C_N =\emptyset$.

Case 2: The wall $\gamma_j^i$ cannot be obtained from weights occuring in the normal form of $\Phi|_{U_l}$ for any such $l$. In this case $q_0$ is a regular point. The openness condition (A) implies that there exists some $q_0'\in Q$ with $\Phi(q_0')\in C_N$. This is again a contradiction.

We conclude that $\Phi(Q)\cap C_b\neq\emptyset$ for all $b=1,\dots,N$. Arguments given before lead to $\overline{C_b}\subseteq \Phi(Q)$ for all $b=1,\dots,N$. But then $B=\overline{C_1}\cup \dots \cup \overline{C_N}\subseteq \Phi(Q)$, and $\Phi(q)$ is indeed an inner point of $\Phi(Q)$.
\end{proof} 

The previous lemma essentially expresses the local equality of $\Phi(M)$ and $\Phi(Q)$. The step towards global equality is now a small one.

\begin{theorem}\label{global}
Let $T\times M \rightarrow M$ be a Hamiltonian torus action on a connected symplectic manifold $M$. Suppose the associated moment map $\Phi:M\rightarrow \t^*$ is proper as a map into a convex open set $\rho \subseteq \t^*$. Let $Q$ be a nonempty closed submanifold of $M$ satisfying two conditions: 
\begin{align*}
(A) & \quad \text{The derivative $D\Phi|_Q(q)$ is surjective at each $q\in Q_{reg}$.} \\
(B) & \quad \text{Let $q\in Q, \ X\in \t_q$. If for each neighborhood $U\subseteq M$ of $q$ there exists} \\
& \quad \text{some $p\in U$ with $\Phi_X(p)>\Phi_X(q)$, then each such $U$ contains some} \\
& \quad \text{$q'\in Q$ with $\Phi_X(q')>\Phi_X(q)$.} 
\end{align*} 
Then, $\ \Phi(Q)=\Phi(M)$.
\end{theorem}
\begin{proof}
Recall that $\Phi(M)$ was shown to be convex in \cite{LMTW}, and we assume that $int(\Phi(M))\neq \emptyset$. In particular, $int(\Phi(M))$ is nonempty and connected. Since $\Phi$ is a closed mapping and $Q$ is closed, it is therefore sufficient to show that $int(\Phi(Q))$ is nonempty, closed and open in $int(\Phi(M))$.
Openness is clear. Since $Q_{reg}\neq \emptyset$, by Lemma \ref{dense}, condition (A) implies $int(\Phi(Q))\neq \emptyset $. To show closedness consider a sequence $\{ \Phi(q_n) \} \subseteq int(\Phi(Q))$ which converges to $\Phi(m)\in int(\Phi(M))$. The image $\Phi(Q)$ is closed, hence $\Phi(m)\in \Phi(Q)\cap int(\Phi(M))$. Then, $\Phi(m)\in int(\Phi(Q))$ follows from Lemma \ref{inner}.
\end{proof}

\begin{remark}
Convexity theorems for $\Phi(M)$ have been established for moment maps $\Phi$ that are not necessarily proper (see e.g. \cite{BOR}). Note that in our proof of Theorem \ref{global} we have made frequent use of the fact that $\Phi$ is a closed mapping. But properness entered mainly to show that a neighborhood of some $\Phi(p)$ is partitioned by finitely many walls $\gamma^i_j$ (see Lemma \ref{ucrit}(2)). Also, convexity of $\Phi(M)$ was only used indirectly. The main ingredient in the proof of Theorem \ref{global} is connectedness of $int(\Phi(M))$. It would be interesting to see under what conditions a non-proper moment map still yields $\Phi(Q)=\Phi(M)$. 
\end{remark}

\section{Applications}\label{appl}

\subsection{Involutions}\label{invol}

The conditions (A) and (B) defined in Section \ref{main} might be hard to check in concrete applications. On the other hand, the submanifold $Q$ often arises as the fixed point set of some involutive diffeomorphism $\sigma:M\rightarrow M$. In this case one might hope to find (natural) conditions on $\sigma$ which ensure that $Q$ satisfies (A) and (B) so that Theorem \ref{global} becomes applicable. Duistermaat's theorem, which was stated in the introduction, provides such an example. In this section we introduce rather weak conditions on $\sigma$ and show that they imply (A) and (B) for the fixed point set $Q$.

\begin{lemma}
Let $T\times M \rightarrow M$ be a Hamiltonian torus action on a connected symplectic manifold $M$ with proper moment map $\Phi:M\rightarrow \rho$. Here, $\rho$ is a convex open subset of $\t^*$. In addition, let $\sigma:M\rightarrow M$ be an involution such that
\begin{align*}
(a) & \quad \Phi \circ \sigma=\Phi \quad \text{and} \\
(b) & \quad \sigma(t.m)=t^{-1}.\sigma(m) \ \forall t\in T, m\in M.
\end{align*} 
Denote by $Q$ the fixed point set of $M$ under $\sigma$ (or a connected component thereof). Then, $Q$ is a closed submanifold of $M$ satisfying conditions (A) and (B) defined in Section \ref{main}.
\end{lemma}
\begin{proof} 
(A) \ Let $q\in Q$ be a $\sigma$-fixed point. Denote by $T_qM=V^+\oplus V^-$ the (+1)-(-1)-eigenspace decomposition of the tangent space with respect to $\sigma$ (we use $\sigma$ to denote the involution on $M$ and its derivative at any point). The invariance property $(a)$ of the moment map implies that for $v\in V^-$,
\[ D\Phi(q).v=(D\Phi(q) \circ \sigma).v=-D\Phi(q).v . \]
Therefore, $D\Phi(q)$ vanishes on $V^-$, and for $u\in V^+, v\in V^-$ one obtains $D\Phi(q).(u+v)=D\Phi(q).u$. In particular, $im(D\Phi(q))=im(D\Phi|_Q(q))$, i.e. $D\Phi|_Q(q)$ is surjective whenever $D\Phi(q)$ is. \\
(B) \ Suppose $q\in Q, \ X\in \t_q$ and for each neighborhood $U\subseteq M$ of $q$ there exists some $p\in U$ with $\Phi_X(p)>\Phi_X(q)$. From the normal form (\ref{nf}) one derives an expression
 for $\Phi_X$ on a neighborhood $U_q$ of $q$,
\[ \Phi_X(v,w,x,y)=\frac{1}{2}\sum_{i=1}^{n} \alpha_i(X) (x_i^2+y_i^2) . \]
Denote all distinct nonzero numbers $\alpha_i(X)$ that occur by $r_1,\dots,r_k$. Due to our general assumption $int(\Phi(M))\neq \emptyset$ there is at least one nonzero $\alpha_i(X)$. The complexified tangent space at $q$ decomposes into eigenspaces with respect to the action of $X$.
\[ (T_qM)_\C=V_0 \oplus \bigoplus_{j=1}^{k} (V^{ir_j}\oplus V^{-ir_j}) . \]
Fix $r_j$ for the moment. Taking derivatives in (b) yields $\sigma(V^{ir_j})\subseteq V^{-ir_j}$. More precisely, a basis $B_+$ for $V^{ir_j}$ is mapped onto a basis $B_-$ for $V^{-ir_j}$. With the appropriate ordering on the basis $B_+\cup B_-$ for $V^{ir_j}+V^{-ir_j}$ the restriction of $\sigma$ to $V^{ir_j}+V^{-ir_j}$ is represented by the matrix
\[ A=\begin{pmatrix} 0 & 1 & & & \\ 1 & 0 & & & \\ & & \ddots & & \\ & & & 0 & 1 \\ & & & 1 & 0 \end{pmatrix} . \]
Changing bases $\sigma|_{V^{ir_j}+V^{-ir_j}}$ takes the matrix expression $B=(\begin{smallmatrix} I_d & 0 \\ 0 & -I_d \end{smallmatrix})$, where $d$ denotes the complex dimension of $V^{ir_j}$. The real matrix $B$ also represents the restriction of $\sigma$ to the real vector space $V_j=T_qM \cap (V^{ir_j}+ V^{-ir_j})$ (note that the Jordan canonical form for the involution $\sigma|_{V_j}$ must be $B$, up to permutation). This shows in particular that the subspace $V_j\cap T_qQ$ of $\sigma$-fixed vectors in $V_j$ is nonempty. In fact, $dim(V_j)=2dim(V_j\cap T_qQ)$. The restrictions of the Hessians of $\Phi_X$ and $\Phi_X|_Q$ to $V_j$ and $V_j\cap T_qQ$, respectively, are given by
\[ H_{\Phi_X}(q)=r_j I_{2d} \quad \text{and} \quad H_{\Phi_X|_Q}(q)=r_j I_{d}. \]
Due to our assumption on points $p\in M$ with $\Phi_X(p)>\Phi(q)$ at least one of the $r_j$ must be positive. Fix a nonzero tangent vector $Z\in V_j \cap T_qQ$ for such a $j$. The curve $\gamma:(-\epsilon,\epsilon)\rightarrow U_q \cap Q$ with $\gamma(0)=q$ and $\gamma'(q)=Z$ contains points $q'\in Q$ that are arbitrarily close to $q$ and satisfy $\Phi(q')>\Phi(q)$.
\end{proof}

\begin{cor}\label{sigmacor}
Let $T\times M \rightarrow M$ be a Hamiltonian torus action on a connected symplectic manifold $M$ with proper moment map $\Phi:M\rightarrow \rho$. Here, $\rho$ is a convex open subset of $\t^*$. In addition, let $\sigma:M\rightarrow M$ be an involution such that
\begin{align*}
(a) & \quad \Phi \circ \sigma=\Phi \quad \text{and} \\
(b) & \quad \sigma(t.m)=t^{-1}.\sigma(m) \ \forall t\in T, m\in M.
\end{align*} 
Denote by $Q$ the fixed point set of $M$ under $\sigma$ (or a connected component thereof). \\
If $Q$ is not empty, then $\Phi(Q)=\Phi(M)$. 
\end{cor}
\begin{proof}
Apply Theorem \ref{global}.
\end{proof}

\begin{remark}
Corollary \ref{sigmacor} shows that the conditions on $\sigma$ in Theorem 3.1 in \cite{KO} (and in Duistermaat's theorem) can be weakened significantly. Compactness of $M$ is not required, and $Q$ need not be Lagrangian. The following Section provides an example in which the full strength of Corollary \ref{sigmacor} is required.  
\end{remark}

\subsection{Van den Ban's theorem}\label{subvdb}

In this section we consider a generalization of Kostant's (nonlinear) convexity theorem which is due to van den Ban \cite{V}. After introducing the necessary notation and stating the theorem we show how it can be formulated in a symplectic framework. The method described here is an extension of the symplectic approach of Lu and Ratiu \cite{LR} towards Kostant's nonlinear convexity theorem. A major difference is that the relevant symplectic manifolds are no longer compact. 

Let $G$ be a real connected semisimple Lie group with finite center, equipped with an involution $\tau$, 
i.e. $\tau$ is a smooth group homomorphism such that $\tau^2=id$. We use the same letter to denote the corresponding involution on the Lie algebra level. We write $H$ for an open subgroup of 
$G^\tau$, the $\tau$-fixed points in $G$. Let $K$ be a $\tau$-stable maximal compact subgroup of $G$. 
The corresponding Cartan involution $\theta$ on $\g$ commutes with $\tau$ and induces the Cartan 
decomposition $\g=\k+\p$. If $\h$ and $\q$ denote the $(+1)$- and $(-1)$-eigenspace of $\g$ 
with respect to $\tau$ one obtains
\[ \g = (\k\cap \h) + (\p\cap \h) + (\k\cap \q) + (\p\cap \q) . \]

We fix a maximal abelian subalgebra $\a^{-\tau}$ of $\p\cap \q$. (In \cite{V} this subalgebra is denoted 
by $\a_{pq}$.) In addition, we choose $\a^{\tau}\subseteq \p\cap \h$ such that $\a := \a^{\tau}+\a^{-\tau}$ 
is maximal abelian in $\p$. Let $\Delta(\g,\a^{-\tau})$ and $\Delta(\g,\a)$ denote the sets of roots for 
the root space decomposition of $\g$ with respect to $\a^{-\tau}$ and $\a$, respectively. We choose a system of positive roots $\Delta^+(\g,\a^{-\tau})$ and a compatible positive system $\Delta^+(\g,\a)$ such that
\[ \Delta^+(\g,\a^{-\tau}) =\{ \alpha |_{\a^{-\tau}} : \alpha\in\Delta^+(\g, \a), \alpha |_{\a^{-\tau}}\ne 0\}. \]
This leads to an Iwasawa decomposition
\[ \g=\n+\a+\k=\n^1+\n^2+\a+\k , \]
where
\begin{eqnarray*}
\n & = & \sum_{\alpha \in \Delta^+(\g,\a)} \g^{\alpha}, \\
\n^1 & = & \sum_{\alpha \in \Delta^+(\g,\a), \alpha |_{\a^{-\tau}}\neq 0} \g^{\alpha} = \sum_{\beta \in \Delta^+(\g,\a^{-\tau})} \g^{\beta}, \\
\n^2 & = & \sum_{\alpha \in \Delta^+(\g,\a), \alpha |_{\a^{-\tau}}= 0} \g^{\alpha}. 
\end{eqnarray*}
Here $\ \g^{\alpha}=\{ X\in \g : [H,X]=\alpha(H)X \ \forall H\in \a \} $ for $\alpha \in \Delta(\g,\a)$, 
and similarly $\g^\beta$ is defined for $\beta \in \Delta(\g,\a^{-\tau})$.  

Let $N$ and $A$ denote the analytic subgroups of $G$ with Lie algebras $\n$ and $\a$, respectively. 
The Iwasawa decomposition $G=NAK$ on the group level has the middle projection $\mu:G\rightarrow A$. 
We write $pr_{\a^{-\tau}}:\a \rightarrow \a^{-\tau}$ for the projection along $\a^{\tau}$. 

For $\beta \in \Delta^+(\g,\a^{-\tau})$ define $H_\beta \in \a^{-\tau}$ such that
\[ H_\beta \perp \mbox{ker} \beta, \quad \beta(H_\beta)=1 , \]
where $\perp$ means orthogonality with respect to the Killing form $\langle , \rangle $. 

Note that the involution $\theta \circ \tau$ leaves each root space 
\[ \g^{\beta}=\sum_{\alpha \in \Delta(\g,\a), \alpha|_{\a^{-\tau}}=\beta} \g^{\alpha} \]
stable. Each $\g^{\beta}=(\g^{\beta})_+ \oplus (\g^{\beta})_- $ decomposes into $(+1)$- and $(-1)$-eigenspace with respect to $\theta \circ \tau$. 

For 
\[ \Delta_{-}:=\{ \beta \in \Delta(\g,\a^{-\tau}) : (\g^{\beta})_- \neq 0 \} , \]
let $\Delta^{+}_{-}=\Delta_{-} \cap \Delta^+(\g,\a^{-\tau})$. Define the closed cone
\[ \Gamma(\Delta^{+}_{-})=\sum_{\beta \in \Delta^{+}_{-}} \R_+ H_\beta . \]
Write $\mathcal{W}_{K \cap H}$ for the Weyl group
\[ \mathcal{W}_{K \cap H}=N_{K \cap H}(\a^{-\tau})/Z_{K \cap H}(\a^{-\tau}) . \]
The convex hull of a Weyl group orbit through $X\in \a^{-\tau}$ will be denoted by $conv(\mathcal{W}_{K \cap H}.X)$. 

We can now state the central theorem.

\begin{theorem}\label{VDB}{\em{(Van den Ban \cite{V})}} \\
Let $G$ be a real connected semisimple Lie group with finite center, equipped with an involution $\tau$, 
and $H$ a connected open subgroup of $G^\tau$. For $X \in \a^{-\tau}$, write $a=\exp X \in A^{-\tau}$. Then
\[ (pr_{\a^{-\tau}} \circ \log \circ \mu) ( H a ) = conv(\mathcal{W}_{K \cap H}.X) + \Gamma(\Delta^{+}_{-}) . \]
\end{theorem}

\begin{remark} \ 
\begin{enumerate}
\item Van den Ban proved his theorem under the weaker condition that $H$ is an essentially connected 
open subgroup of $G^\tau$ (by reducing it to the connected case). 
\item If $\tau=\theta$ one obtains Kostant's (nonlinear) convexity theorem. Note that in this case 
the group $H$ and the orbit $Ha$ are compact.
\end{enumerate}
\end{remark}

The map $pr_{\a^{-\tau}} \circ \log \circ \mu:Ha \longrightarrow \a^{-\tau}$ clearly factors to a map on the orbit $Ha/K$ in the symmetric space $G/K$. In this setting van den Ban's theorem consists of a precise description of the image $(pr_{\a^{-\tau}} \circ \log \circ \mu)(Ha/K)$ in $\a^{-\tau}$.

In the special case $G$ complex and $\tau$ complex antilinear the orbit $Ha/K \subseteq G/K$ carries a Poisson structure which was introduced in a much more general context by Evens and Lu \cite{EL}. The symplectic leaf in $Ha/K$ through $aK$ turns out to be diffeomorphic to the coset space $A^\tau \setminus Ha/K$ (cf. \cite{FO}). Moreover, the action of the torus $\exp(i\a^{-\tau})$ on $A^\tau \setminus Ha/K$ induced from multiplication from the left is Hamiltonian with moment map $pr_{\a^{-\tau}} \circ \log \circ \mu:A^\tau \setminus Ha/K \longrightarrow \a^{-\tau}$. (Note that $A^\tau=\exp (\a^\tau)$ normalizes $\n$ if $G$ is complex.) This moment map is proper as was already shown in \cite{V}. A version of the Atiyah-Guillemin-Sternberg convexity theorem for noncompact symplectic manifolds \cite{HNP} now implies that the moment map image $(pr_{\a^{-\tau}} \circ \log \circ \mu)(A^\tau \setminus Ha/K)=(pr_{\a^{-\tau}} \circ \log \circ \mu)(HaK)$ is convex and locally polyhedral. A careful analysis of the local moment cones shows that this image is in fact the sum of the compact convex set and the cone in van den Ban's theorem (see \cite{FO}).

\begin{remark}
In case that $\tau=\theta$, the subgroup $A^\tau$ is trivial. Therefore, $Ha/K=Ka/K$ is symplectic. Its symplectic structure coincides with the one induced by the Lu-Weinstein Poisson structure on $G/K$. 
\end{remark}

For real $G$ the symplectic approach just described is not directly applicable. It is no loss of generality (see \cite{V}, Lemma 2.1) to assume that $G$ admits a (simply connected) complexification $G_\C$ with Lie algebra $\g_\C$. Now $\tau$ and $\theta$ will denote the involutions on $\g$ and $G$, as well as their antilinear extensions to $\g_\C$ and $G_\C$. To avoid ambiguity in the notation we will frequently add tildes when we refer to the corresponding object in the setting of $G_\C$. For example, the Cartan decomposition of $\g_\C$ becomes $\g_\C=\tilde{\k}+\tilde{\p}$, where $\tilde{\k}=\k+i\p$ and $\tilde{\p}=i\k+\p$. The decomposition with respect to $\tau$ becomes $\g_\C=\tilde{\h}+\tilde{\q}$.

Let $\a^{-\tau}$ be the maximal abelian subalgebra of $\p^{-\tau}$. We fix $\r \subseteq \k^\tau$ such that $\tilde{\a}^{-\tau}=\a^{-\tau}+i\r$ is maximal abelian in $\tilde{\p}^{-\tau}$. Also we can find $\tilde{\a}^\tau $ such that $\tilde{\a}=\tilde{\a}^{-\tau}+\tilde{\a}^\tau $ is maximal abelian in $\tilde{\p}$. We choose a positive system $\Delta^+(\g,\tilde{\a})$ related to an ordering of a basis for $\tilde{\a}$, where a basis for $\a^{-\tau}$ comes first, followed by a basis for $i\r$, and finally a basis for $\tilde{\a}^\tau$. The relevant Iwasawa decompositions $G=NAK$ and $G_\C=\tilde{N}\tilde{A}\tilde{K}$ have middle projections $\mu$ and $\tilde{\mu}$. Note that in general $A \not\subseteq \tilde{A}$. To compare $\mu$ and $\tilde{\mu}$ we introduce the projection $pr:\tilde{\a}^{-\tau}\rightarrow \a^{-\tau}$ along $i\r$.

\begin{lemma}\label{mus}
For each $x\in G$, 
\[ pr_{\a^{-\tau}} \circ \log \circ \mu(x)=pr \circ pr_{\tilde{\a}^{-\tau}} \circ \log \circ \tilde{\mu}(x) . \]
\end{lemma}
\begin{proof}
Recall
\begin{eqnarray*}
\n^1 & = & \sum_{\alpha \in \Delta^+(\g,\a), \alpha|_{\a^{-\tau}}\neq 0} \g^\alpha \subseteq \n , \\
\n^2 & = & \sum_{\alpha \in \Delta^+(\g,\a), \alpha|_{\a^{-\tau}}=0} \g^\alpha \subseteq \n , 
\end{eqnarray*}
and $N^1$ and $N^2$, the analytic subgroups of $G$ with Lie algebras $\n^1$ and $\n^2$, respectively. \\
For $x\in G$ we can write
\[ x=mnabk \qquad {\mbox with} \quad m\in N^1, n\in N^2, a\in A^\tau, b\in A^{-\tau}, k\in K . \]
We want to show that $\ \log b=pr \circ pr_{\tilde{\a}^{-\tau}} \circ \log \circ \tilde{\mu}(x)$. 

Since $N^1\subseteq \tilde{N}$, due to our choice of positive systems, and $K\subseteq \tilde{K}$, we can assume that $m=k=e$. But $x=nab$ lies in $Z(\a^{-\tau})$, the analytic subgroup of $G_\C$ with Lie algebra $\z_{\g_\C}(\a^{-\tau})$, which is reductive. Write 
$\ \z_{\g_\C}(\a^{-\tau})=\g_{ss}+\z \ $, with $\g_{ss}=[\z_{\g_\C}(\a^{-\tau}), \z_{\g_\C}(\a^{-\tau})]$ and $\z$ the center of $\z_{\g_\C}(\a^{-\tau})$. Let $G_{ss}$ and $Z$ denote the corresponding analytic subgroups of $G_\C$. 

We consider the hermitian inner product $B(X,Y)=-\Re \kappa(X,\theta Y)$ on $\g_\C$ and possible restrictions to ($\theta$-stable) subspaces. One verifies $\a^\tau \perp_B \a^{-\tau}$ and $\g_{ss} \perp_B \z$. Then,
\[ \z_{\g_\C}(\a^{-\tau})=\g_{ss}+(\z \cap (\a^{-\tau})^\perp)+\a^{-\tau} . \]
The Lie algebra $\l=\g_{ss}+(\z \cap (\a^{-\tau})^\perp)$ is reductive, $\theta$-stable, and it contains $\n^2+\a^\tau$ and $\tilde{a}\cap (\a^{-\tau})^\perp$. Moreover, $\tilde{a}\cap (\a^{-\tau})^\perp$ is maximal abelian in $\l \cap \tilde{p}$. The analytic subgroup $L\subseteq G_\C$ with Lie algebra $\l$ is reductive and admits an Iwasawa decomposition
\[ L=N_L A_L K_L \quad {\mbox with} \ N_L\subseteq \tilde{N}, \ A_L=\exp(\tilde{a}\cap (\a^{-\tau})^\perp)=\exp(\tilde{\a}^\tau+i\r) \subseteq \tilde{A}, \ K_L \subseteq \tilde{K} . \]
Accordingly write $na\in L$ as $na=n_La_Lk_L$ with $n_L \in N_L, a_L \in A_L, k_L \in K_L$. Then
\[ \tilde{\mu}(x)=\tilde{\mu}(nab)=\tilde{\mu}(n_La_Lbk_L)=a_L b , \]
and therefore
\[ pr \circ pr_{\tilde{\a}^{-\tau}} \circ \log \circ \tilde{\mu}(x)=\log b . \]
\end{proof}

We fix $X\in \a^{-\tau}$ and $a=\exp X \in A^{-\tau}$. Van den Ban's theorem for the pair $(G_\C,\tau)$ states that
\begin{equation}\label{vdbceq} 
(pr_{\tilde{\a}^{-\tau}} \circ \log \circ \tilde{\mu}) ( \tilde{H} a ) = conv(\mathcal{W}_{\tilde{K} \cap \tilde{H}}.X) + \Gamma(\tilde{\Delta}^{+}_{-}) , 
\end{equation}
where $ \tilde{\Delta}^+_-= \{ \beta \in \Delta^{+} (\g_\C ,\tilde{\a}^{-\tau}) : (\g_\C^\beta)_- \neq 0 \} $. (Note that this result can be obtained by ''symplectic means'' \cite{FO} as explained above.)

\begin{lemma}\label{cones}
\[ pr(conv(\mathcal{W}_{\tilde{K} \cap \tilde{H}}.X) + \Gamma(\tilde{\Delta}^{+}_{-}))=conv(\mathcal{W}_{K \cap H}.X) + \Gamma(\Delta^{+}_{-}) . \]
\end{lemma}
\begin{proof}
We consider the compact and the conal part separately. 

To see that $pr(conv(\mathcal{W}_{\tilde{K} \cap \tilde{H}}.X))=conv(\mathcal{W}_{K \cap H}.X)$ note that this is a statement about the reductive Lie algebra $\g^{\theta \tau}\subseteq \g$ of $\theta \tau$-fixed points and its complexification. The restrictions of $\theta$ and $\tau$ to $\g^{\theta \tau}$ coincide, and the Cartan decomposition becomes $\g^{\theta \tau}=\k^\tau + \p^{-\tau}$. In particular, $\a^{-\tau}$ is maximal abelian in $\p^{-\tau}$. Without loss of generality we can assume that $\g^{\theta \tau}$ is semisimple. Complex conjugation $Z\mapsto \bar{Z}$ on $\g^{\theta \tau}_\C$ with respect to the real form $\g^{\theta \tau}$ induces an involution on $\tilde{\a}^{-\tau}$ with fixed point set $\a^{-\tau}$. The map $pr:\tilde{\a}^{-\tau}\rightarrow \a^{-\tau}$ then just returns the real part of its argument. Since $conv(\mathcal{W}_{\tilde{K} \cap \tilde{H}}.X)$ is convex and invariant under this involution, for $x\in conv(\mathcal{W}_{\tilde{K} \cap \tilde{H}}.X)$, we obtain 
\[ pr(x)=\frac{x+\bar{x}}{2}\in conv(\mathcal{W}_{\tilde{K} \cap \tilde{H}}.X)\cap \a^{-\tau} \subseteq conv(\mathcal{W}_{K \cap H}.X) . \]
This proves $pr(conv(\mathcal{W}_{\tilde{K} \cap \tilde{H}}.X))\subseteq conv(\mathcal{W}_{K \cap H}.X)$. The other inclusion is trivial. 

To show $pr(\Gamma(\tilde{\Delta}^{+}_{-}))=\Gamma(\Delta^{+}_{-})$, first consider $\beta\in \tilde{\Delta}^{+}_{-}$, i.e. $\beta \in \Delta^+(\g_\C,\tilde{\a}^{-\tau})$ with $\theta \tau Y_\beta = -Y_\beta $ for some $Y_\beta \neq 0$. We can write $Y_\beta=A_\beta +iB_\beta$ with $A_\beta, B_\beta \in \g$. Since $\beta$ is real on $\a^{-\tau}$ and since both $\theta$ and $\tau$ are antilinear there is some $C_\beta \in \g$ (choose $A_\beta$ or $B_\beta$, whichever is nonzero) with 
\[ \theta \tau C_\beta = -C_\beta , \quad [Z,C_\beta]=\beta(Z)C_\beta \quad \forall \ Z\in \a^{-\tau} . \]
But this is just the condition that $\beta |_{\a^{-\tau}} \in \Delta^+_- $. In other words, $pr(H_\beta)=H_{\beta |_{\a^{-\tau}}}\in \Gamma(\Delta^{+}_{-})$. This shows $pr(\Gamma(\tilde{\Delta}^{+}_{-}))\subseteq \Gamma(\Delta^{+}_{-})$. 
Conversely, if $\gamma \in \Delta^+_- $ then there is $C_\gamma \in \g$ with 
\[ \theta \tau C_\gamma = -C_\gamma , \quad [Z,C_\gamma]=\gamma(Z)C_\gamma \quad \forall \ Z\in \a^{-\tau} . \]
The restricted root space $\g_\C^\gamma \subseteq \g_\C$ decomposes as $\g_\C^\gamma = \sum_{\beta \in \Delta(\g_\C,\tilde{\a}^{-\tau}), \beta |_{\a^{-\tau}}=\gamma} \g_\C^\beta $, therefore we can write $C_\gamma=\sum_{\beta_i |_{\a^{-\tau}}=\gamma} C_{\beta_i}$ with $C_{\beta_i}\in \g_\C^{\beta_i}$. Since each $\g_\C^\beta$ is $\theta \tau$-stable there is some $\beta_i$ with $\beta_i |_{\a^{-\tau}}=\gamma$ and a nonzero $C_{\beta_i}\in \g_\C^{\beta_i}$ with $\theta \tau C_{\beta_i}=-C_{\beta_i}$. But this means $\beta_i \in \tilde{\Delta}^{+}_{-}$. Hence $H_\gamma=H_{\beta_i |_{\a^{-\tau}}}=pr(H_{\beta_i})\in pr(\Gamma(\tilde{\Delta}^{+}_{-}))$, which implies $\Gamma(\Delta^{+}_{-}) \subseteq pr(\Gamma(\tilde{\Delta}^{+}_{-}))$.
\end{proof}

We can now finish the symplectic proof of van den Ban's theorem. 

\begin{theorem}\label{vdbr}(van den Ban) 
\[ pr_{\a^{-\tau}} \circ \log \circ \mu(Ha)=conv(\mathcal{W}_{K \cap H}.X) + \Gamma(\Delta^{+}_{-}) . \]
\end{theorem}
\begin{proof}
In view of Lemmas \ref{mus} and \ref{cones} and the complex version of van den Ban's theorem (\ref{vdbceq}) it is sufficient to show
\begin{equation}\label{Qeq}
pr \circ pr_{\tilde{\a}^{-\tau}} \circ \log \circ \tilde{\mu}(Ha)=pr \circ pr_{\tilde{\a}^{-\tau}} \circ \log \circ \tilde{\mu}(\tilde{H}a) .
\end{equation}
Recall that the right hand side equals the image of the moment map $\Phi=pr \circ pr_{\tilde{\a}^{-\tau}} \circ \log \circ \tilde{\mu}$ for the Hamiltonian action of $T=\exp(i\a^{-\tau})$ on the symplectic manifold $M=\tilde{A}^\tau\setminus \tilde{H}a/\tilde{K}$. Complex conjugation $Z\mapsto \bar{Z}$ on $\g_\C$ with respect to $\g$ leaves $\tilde{\a}^\tau$ and $\tilde{\k}$ invariant and therefore induces an involution $\sigma$ on $M$. We denote the connected component of the fixed point set $M^\sigma$ through the point $\tilde{A}^\tau a\tilde{K}\in M$ by $Q$. Clearly, $Q=\tilde{A}^\tau\setminus Ha/\tilde{K}$ is closed and nonempty. Now, equation (\ref{Qeq}) is equivalent to $\Phi(Q)=\Phi(M)$. In order to apply Corollary \ref{sigmacor} we need to show that $\sigma:M\rightarrow M$ satisfies two properties.

(a) \ $\Phi \circ \sigma =\Phi $. \\
We proceed in a similar fashion as in the proof of Lemma \ref{mus}. Define the following subalgebras of $\g_\C$.
\begin{eqnarray*}
\tilde{\n}^1 & = & \sum_{\alpha \in \Delta^+(\g_\C,\tilde{\a}), \alpha|_{\a^{-\tau}}\neq 0} \g_\C^\alpha \subseteq \tilde{\n} , \\
\tilde{\n}^2 & = & \sum_{\alpha \in \Delta^+(\g_\C,\tilde{\a}), \alpha|_{\a^{-\tau}}=0} \g_\C^\alpha \subseteq \tilde{\n} , 
\end{eqnarray*}
with corresponding analytic subgroups $\tilde{N}^1, \tilde{N}^2$ of $G_\C$. 
Then $x\in G_\C$ can be written as
\[ x=mnabk \quad \mbox{with} \ m\in \tilde{N}^1, n\in \tilde{N}^2, a\in \exp(\tilde{\a}^\tau+i\r), b\in A^{-\tau}, k\in \tilde{K} . \]
Then $\ \sigma(x)=\sigma(m)\sigma(na)b\sigma(k)$. Note that $\sigma(m)\in \tilde{N}^1\subseteq \tilde{N}, \sigma(k)\in \tilde{K}$. On the Lie algebra level $\sigma$ leaves $\z_{\g_\C}(\a^{-\tau})$ stable. Since $\a^{-\tau}$ is $\sigma$-stable and $B$ is $\sigma$-invariant, $\sigma$ also leaves $\l$ stable (recall the notation from Lemma \ref{mus}). Therefore, $\sigma(na)\in L$ decomposes as 
\[ \sigma(na)=n_L a_L k_L \qquad {\mbox with} \quad n_L \in N_L, a_L\in A_L=\exp(\tilde{\a}^\tau+i\r), k_L \in K_L . \]
This means $\sigma(x)=\sigma(m)n_L a_L b k_L \sigma(k)$, i.e.
\[ \Phi \circ \sigma (x)=pr \circ pr_{\tilde{\a}^{-\tau}} \circ \log \circ \tilde{\mu}(x)=\log b=\Phi(x) . \]

(b) \ $t\circ \sigma = \sigma \circ t^{-1} \quad \forall \ t\in T . $ \\
Fix $h\in \tilde{H}$ and $t\in T=\exp(i\a^{-\tau})$. Then,
\[ (\sigma \circ t^{-1})(\tilde{A}^\tau ha\tilde{K})=\tilde{A}^\tau \overline{t^{-1}ha}\tilde{K})=\tilde{A}^\tau t\bar{h}a\tilde{K}=t \circ \sigma(\tilde{A}^\tau ha\tilde{K}) . \]

Now that properties (a) and (b) of $\sigma$ have been established Corollary \ref{sigmacor} can be applied. It asserts $\Phi(Q)=\Phi(M)$ or
\[ pr \circ pr_{\tilde{\a}^{-\tau}} \circ \log \circ \tilde{\mu}(\tilde{A}^\tau\setminus Ha/ \tilde{K})=pr \circ pr_{\tilde{\a}^{-\tau}} \circ \log \circ \tilde{\mu}(\tilde{A}^\tau\setminus \tilde{H}a/ \tilde{K}) . \]
This implies (\ref{Qeq}) and concludes the proof.

\end{proof}

\bibliographystyle{amsplain}

\end{document}